\documentclass[11pt]{article}
\usepackage{amsmath}
\usepackage{amssymb}
\usepackage{epsfig}
\def\BC{{\mathbb C}}
\def\BD{{\mathbb D}}

\newcommand{\bpr}{{\noindent\textbf{Proof.}\ }}
\newcommand{\epr}{{\hfill $\Box$}}
\newcommand{\K}{K}

\newcommand{\sH}{{\cal H}}
\newcommand{\sK}{{\cal K}}

\newcommand{\sF}{{\cal F}}
\newcommand{\sU}{{\cal U}}
\newcommand{\sS}{{\cal S}}

\newcommand{\sG}{{\cal G}}

\newcommand{\sE}{{\cal E}}
\newcommand{\sD}{{\cal D}}

\newcommand{\sY}{{\cal Y}}

\newcommand{\im}{\textup{Im\,}}
\newcommand{\mat}[2]{\ensuremath{\left[\begin{array}{#1}
#2
\end{array} \right]}}

\newcommand{\la}{\lambda}
\newcommand{\Om}{\Omega}
\newcommand{\om}{\omega}

\newcommand{\eL}{{\mathbf{L}}}
\newcommand{\eS}{{\mathbf{S}}}
\newcommand{\eD}{{\mathbf{D}}}

\newcommand{\col}{\textup{col\,}}

\newcommand{\opendisc}{\mathbb{D}}
\newcommand{\closedisc}{\overline{\mathbb{D}}}

\newcommand{\half}{\frac{1}{2}}

\newtheorem{theorem}{Theorem}[section]

\newtheorem{lemma}[theorem]{Lemma}
\newtheorem{proposition}[theorem]{Proposition}

\begin{document}

\title{All solutions to the relaxed commutant\\ lifting problem}
\author{A.E. Frazho\footnote{
The research of the first author was supported in part by a visitor's grant
from NWO (Nederlandse Organisatie voor Wetenschappelijk Onderzoek).},
S. ter Horst and M.A. Kaashoek}
\date{}
\maketitle

\begin{abstract}
A new description is given of  all solutions to the relaxed
commutant lifting problem. The method of proof is also different
from earlier ones, and uses only an operator-valued version of a
classical lemma on harmonic  majorants.
\end{abstract}



\setcounter{section}{-1}
\section{Introduction}
In this paper we give a new, more refined, description of  all solutions to
the relaxed commutant lifting problem. Let us first recall the
formulation of this problem. The starting point is a \emph{data
set} $\{A,T^\prime,U^\prime,R,Q\}$ consisting of five Hilbert
space operators. The operator $A$ is a contraction mapping $\sH$
into ${\sH}'$, the operator $U^\prime$ on $\sK^\prime$ is a
minimal isometric lifting of the contraction $T^\prime$ on $\sH'$,
and $R$ and $Q$ are operators from $\sH_0$ to $\sH$, satisfying the
following constraints:
\begin{equation}\label{intertw}
T'AR=AQ  \quad \textup{and}\quad R^*R\leq Q^*Q.
\end{equation}
Given this data set the relaxed commutant lifting problem
(\emph{RCL problem}) is to find all contractions $B$ from $\sH$ to
$\sK'$ such that
\begin{equation}\label{rclt}
\Pi_{\sH'}B=A\quad\mbox{and}\quad U'BR=BQ.
\end{equation}
Here $\Pi_{\sH'}$ is the orthogonal projection from $\sK'$ onto
$\sH'$.

The RCL problem has been introduced in \cite{ffk02}, and in the
paper \cite{ffk02} also an explicit construction for a particular
solution is given. By choosing $\sH_0 = \sH$ with
$R$ the identity operator on $\sH$ and $Q = T$ an isometry
on $\sH$, one sees that the solution of the RCL
problem in \cite{ffk02} contains the classical Sz-Nagy-Foias
commutant lifting theorem \cite{Sz.-NF68} as a special case. Also
a number of recent generalizations of the commutant lifting
theorem can be seen as special cases of the solution to the RCL
problem presented in \cite{ffk02}. This includes the Treil-Volberg
version \cite{TV94}, which appears when one takes $R = I$, and the
weighted commutant lifting theorem from \cite{BFF99}. Finally,
\cite{ffk02} also shows that the solution of the RCL problem allows
one to solve relaxed versions of most metric constrained
interpolation problems from \cite{FFGK98}, and  their $H^2$ versions.

In \cite{FtHK04} a Redheffer type description is given of all
solutions to the RCL problem by using the theory of isometric
realizations and Arocena's coupling method from
\cite{ar83a} and \cite{ar83e}, see also Section VII.8 in \cite{ff}. A
choice sequence approach for the description of all solutions,
also using the coupling framework, can be found in \cite{LT04}. In
the present paper we give a more refined and more explicit
description of all solutions than the one appearing in
\cite{FtHK04}. Furthermore, our proof will be rather elementary
and uses only an operator-valued version of a classical result on
harmonic majorants. Our approach is even interesting in the
classical commutant lifting setting, and provides a new proof for
Theorem XIII.3.4 in \cite{ff} (see the final part of Section
\ref{sec:mainth}).

The paper consists of three sections not counting the present
introduction. In the first section we introduce the necessary
terminology, state our two main theorems, and specify our results
for the commutant lifting setting. The second section contains
preliminary material on positive real operator-valued functions
and presents an operator-valued version of a classical result on
least harmonic majorants (cf., \cite{D70}, page 28). The proofs of
our two main theorems are given in the third section.

We conclude with a few words about notation and terminology.
Throughout capital calligraphic letters denote Hilbert spaces. The
Hilbert space direct sum of $\sU$ and $\sY$ is denoted by
\[
\sU\oplus \sY \quad \mbox{or by}\quad \left[
\begin{array}{c}
\sU\\ \sY
\end{array}
\right].
\]
The term \emph{operator} stands for a bounded linear
transformation acting between Hilbert spaces. The set of all
operators from $\sU$ into $\sY$ is denoted by $\eL(\sU, \sY)$. The
identity operator on the space $\sU$ is denoted by $I_{\sU}$ or
just by $I$, when the underlying space is clear from the context.
As usual, given a contraction $A$ from $\sU$ into $\sY$, we write
$D_A$ for the defect operator $(I_\sU-A^*A)^{1/2}$ and $\sD_A$ for
the closure of the range of $D_A$. For the definition of an
isometric lifting and a review of its properties we refer to
Section IV.1 in \cite{FFGK98}. By definition, a $\eL(\sU,
\sY)$-valued \emph{Schur class function} is a function which is
analytic on the open unit disk $\BD$ and whose values are
contractions from $\sU$ to $\sY$. The class of these functions is
denoted by $\eS(\sU,\sY)$ and is called a \emph{Schur class}.
Notice that a function $F$ belongs to the Schur class
$\eS(\sE,\sY_1\oplus\sY_2)$ if and only if $F$ admits a matrix
representation of the form
\begin{equation}\label{Frep}
F(\la)=\mat{c}{F_1(\la)\\F_2(\la)},\quad \la\in\BD,
\end{equation}
where $F_1$ is in $\eS(\sE,\sY_1)$ and $F_2$ is in
$\eS(\sE,\sY_2)$ such that
$F_1(\la)^*F_1(\la)+F_2(\la)^*F_2(\la)\leq I$ for all $\la\in\BD$.
For convenience a function $F$ that is represented as in
(\ref{Frep}) will be denoted by $F=\col[F_1,F_2]$. By $H^2(\sU)$
we denote the Hardy space of all $\sU$-valued analytic functions
$f$ on $\mathbb{D}$ such that $\sum_{\nu=0}^\infty
\|f_\nu\|^2<\infty$, where $f_0,f_1,f_2, \dots$ are the Taylor
coefficients of $f$ at zero. Finally, $S_\sU$ denotes the
unilateral shift on $H^2(\sU)$ and $E_\sU$ is the canonical
embedding of $\sU$ onto the space of constant functions in
$H^2(\sU)$ defined by $(E_\sU v)(\lambda)\equiv v$ for all $v\in
\sU$. We simply write $S$ and $E$ if the underlying space is clear
from the context.

\section{Main theorems}\label{sec:mainth}
\setcounter{equation}{0} Let $\{A,T^\prime,U^\prime,R,Q\}$ be a
fixed data set. In the sequel we say that $B$ is a \emph{solution
to the RCL problem} for the data set $\{A,T^\prime,U^\prime,R,Q\}$
if $B$ is a contraction from $\sH$ into $\sK'$ satisfying
(\ref{rclt}).

Without loss of generality we shall assume that $U^\prime$ is the
Sz.-Nagy-Sch\"{a}ffer minimal isometric lifting of $T^\prime$,
that is,
\begin{equation}\label{szns}
 U^\prime = \left[\begin{array}{cc}
   T^\prime & 0 \\
   ED_{T^\prime} & S \\
 \end{array}\right] \mbox{ on } \sK' =  \left[\begin{array}{c}
   \sH^\prime \\
   H^2(\sD_{T^\prime}) \\
 \end{array}\right].
\end{equation}
Here $S$ is the unilateral shift on $H^2(\sD_{T^\prime})$ and $E$
is the canonical embedding of $\sD_{T^\prime}$ into
$H^2(\sD_{T^\prime})$ defined by $(Ed)(\lambda)\equiv d$ for all
$d\in \sD_{T^\prime}$.

Since we assume that $\sK'=\sH'\oplus H^2(\sD_{T'})$, an operator $B$ from
$\sH$ into $\sK'$ is a contraction satisfying $\Pi_{\sH'}B=A$, as in the first
identity of~(\ref{rclt}), if and only if $B$ can be represented in the form
\begin{equation}\label{defB}
B= \left[\begin{array}{c}
     A  \\
   \Gamma D_A \\
 \end{array}\right]: \sH \to \left[\begin{array}{c}
   \sH^\prime \\
   H^2(\sD_{T^\prime}) \\
 \end{array}\right],
\end{equation}
where $\Gamma$ is a contraction from $\sD_A$ into $H^2(\sD_{T^\prime})$.
Moreover, $B$ and $\Gamma$ determine each other uniquely.
Using this representation of $B$ and the fact that $U'$ is given
by~(\ref{szns}), the constraint $U'BR=BQ$ in (\ref{rclt}) is equivalent to
\begin{equation}\label{fundeq}
E D_{T^\prime}AR+S\Gamma D_AR=\Gamma D_AQ.
\end{equation}
Therefore, with $U'$ as in (\ref{szns}), the RCL problem for
$\{A,T^\prime,U^\prime,R,Q\}$ is equivalent to the problem of
finding all contractions $\Gamma$ from $\sD_A$ into
$H^2(\sD_{T^\prime})$ such that (\ref{fundeq}) holds.

To state our two main theorems we need some additional notation.
Observe that, because of (\ref{intertw}), for each $h\in \sH_0$ we
have
\begin{eqnarray}\label{om3}
\|D_A Q h\|^2 &=& \|Q h\|^2 - \|A Q h\|^2 \geq
\|R h\|^2 - \|T^\prime A R h\|^2\nonumber\\
&=&  \|A R h\|^2 - \|T^\prime A R h\|^2 + \|R h\|^2 -\|A R h\|^2\nonumber\\
&=& \|D_{T^\prime} A R h\|^2 + \|D_A R h\|^2.
\end{eqnarray}
Hence the identity
 \begin{equation}\label{535}
\omega D_A Q h = \left[\begin{array}{c}
  D_{T^\prime} A R h \\
  D_A R h\\
\end{array}\right],\quad h \in \sH_0,
\end{equation}
uniquely defines a contraction $\omega$ from $\sF=\overline{D_A Q
\sH}$ into $\sD_{T^\prime} \oplus \sD_{A}$. Let  $\omega_1$ be the
contraction mapping $\sF$ into $\sD_{T^\prime}$
 determined by the first component of $\omega$
and $\omega_2$ be the contraction mapping $\sF$ into $\sD_{A}$
determined by the second component of $\omega$, that is,
\begin{equation}\nonumber
\omega_1 D_A Q h =  D_{T^\prime} A R h \quad \mbox{and}\quad
\omega_2 D_A Q h = D_A R h, \quad\mbox{for all } h \in \sH_0.
\end{equation}
Notice that we have equality in (\ref{om3}) if and only if  $R^*R
= Q^*Q$. In other words, $\omega$ is an isometry if and only if
$R^*R = Q^*Q$, which happens in many applications. In particular,
$\omega$ is an isometry in the setting of the commutant lifting
problem.
The equation in (\ref{fundeq}) can equivalently be represented in terms of
$\om_1$ and $\om_2$ as
\begin{equation}\label{altfundeq}
E\om_1+S\Gamma\om_2=\Gamma|\sF.
\end{equation}

In the sequel we shall call a pair of operator-valued functions
$\{F,G\}$ a \emph{Schur pair associated with the data set}
$\{A,T^\prime,U^\prime,R,Q\}$ if $\col[F,G]$ is in
$\eS(\sD_A,\sD_{T'}\oplus\sD_A)$ and $\col[F,G](\la)|\sF=\om$ for
all $\la\in\BD$. In other words, $\{F,G\}$ is a Schur pair if both
$F$ and $G$ are analytic operator-valued functions, where
$F:\mathbb{D}\to \eL(\sD_A,\sD_{T^\prime})$ and $G:\mathbb{D}\to
\eL(\sD_A,\sD_A)$, such that
\begin{equation}\label{schurpair}
F(\lambda)^*F(\lambda)+G(\lambda)^*G(\lambda)\leq I,\quad
F(\lambda)|\sF = \omega_1, \quad
 G(\lambda)|\sF = \omega_2,
 \quad\mbox{for all } \lambda \in \mathbb{D}.
 \end{equation}
We can now state the first main theorem.

\begin{theorem} \label{thrc}
Consider the data set $\{A,T^\prime,U^\prime,R,Q\}$ with
$U^\prime$ being given by $(\ref{szns})$. Then all solutions to the
corresponding RCL problem are given by
\begin{equation}\label{sols1}
B=\mat{c}{A\\\Gamma D_A}:\sH\to\mat{c}{\sH'\\H^2(\sD_{T'})},
\end{equation}
where $\Gamma$ is a contraction from $\sD_A$ into $H^2(\sD_{T'})$ given by
\begin{equation}\label{sols2}
(\Gamma d)(\lambda)=F(\lambda)(I-\lambda G(\lambda))^{-1}d,\quad
d\in\sD_A,\lambda\in\BD,
\end{equation}
with $\{F,G\}$ an arbitrary Schur pair associated with the given
data set.
\end{theorem}

The mapping $\{F,G\}\mapsto B$ from the set of Schur pairs to the
solutions of the RCL problem described in Theorem \ref{thrc} is
onto but not necessarily one to one. In other words, in general
there can by  many Schur pairs associated with a specified
solution $B$, via (\ref{sols1}) and (\ref{sols2}). However, in the
classical commutant lifting setting the mapping $\{F,G\}\mapsto B$
is onto and one to one, see \cite{ff} and the final part of this
section. To describe the non-uniqueness we need some additional
notation.

Let $B$ in (\ref{sols1}) be a fixed solution to the RCL problem for the
data set $\{A,T^\prime,U^\prime,R,Q\}$ with $U^\prime$
being given by $(\ref{szns})$, and let $\Gamma$ be the
contraction from  $\sD_A$ into $H^2(\sD_{T^\prime})$ determined by
$B$ via (\ref{sols1}). Then $\Gamma$ satisfies (\ref{fundeq}). This
implies that
 there exists a contraction
 $\Omega$ mapping $\sF_\Gamma = \overline{D_\Gamma \sF}$
 into $\sD_\Gamma$ satisfying
 \begin{equation}\label{omg}
\Omega D_\Gamma D_A Q h = D_\Gamma D_A R h, \quad h \in \sH_0.
\end{equation}
To see this we use (\ref{fundeq}) and (\ref{535}) to show that for
all $h$ in $\sH_0$, we have
\begin{eqnarray}\label{739}
\|D_\Gamma D_A Q h\|^2 &=& \|D_A Qh\|^2 - \|\Gamma D_AQ  h\|^2
= \|D_A Qh\|^2 - \|E D_{T'}AR h\|^2 - \|S\Gamma D_AR h\|^2\nonumber\\
&=& \|D_A Qh\|^2 - \|D_{T'}AR h\|^2 - \|\Gamma D_AR h\|^2\nonumber\\
&=&  \|D_{\Gamma}D_A R h\|^2 + \|D_A Qh\|^2 - \|D_{T'}AR h\|^2 -\|D_A R h\|^2\nonumber\\
&=&  \|D_{\Gamma}D_A R h\|^2 + \|D_A Qh\|^2 - \|\omega D_A Q
h\|^2\nonumber\\
&=&  \|D_{\Gamma}D_A R h\|^2 + \|D_\omega D_A Qh\|^2.\nonumber\\
&\geq &\|D_{\Gamma}D_A R h\|^2.
\end{eqnarray}
Thus $\|D_\Gamma D_A Q h\| \geq  \|D_\Gamma D_A R h\|$ for all $h$
in $\sH_0$. So the relation $\Omega D_\Gamma D_A Q  = D_\Gamma D_A
R$ uniquely defines a contraction from $\sF_\Gamma =
\overline{D_\Gamma \sF}$  into $\sD_\Gamma$. By employing the definition
of $\om$ observe that for all $f\in\sF$ we have
$\Omega D_\Gamma f=D_\Gamma\om_2 f$. {}From the calculation
leading to (\ref{739}) we also see that $\Omega$ is an isometry if
and only if $\omega$ is an isometry, and as we saw the latter happens if
and only if $R^*R = Q^*Q$. In particular, $\Omega$ is an isometry
in the setting of the commutant lifting theorem.

Now for $\Gamma$ and $\Omega$ as in the previous paragraph, let
$\eS_\Omega(\sD_\Gamma,\sD_\Gamma)$ be the subset of the Schur
class $\eS(\sD_\Gamma,\sD_\Gamma)$ defined by
\begin{equation}\label{sxi}
\eS_\Omega(\sD_\Gamma,\sD_\Gamma) = \left\{C \in
\eS(\sD_\Gamma,\sD_\Gamma)\, :\,  C(\lambda)|\sF_\Gamma = \Omega
\mbox{ for all } \lambda \in \mathbb{D}\right\}.
\end{equation}
Notice that $\eS_\Omega(\sD_\Gamma,\sD_\Gamma)$ is not empty.  For
example, it contains the function $C$ given by
$C(\la)=\Omega\Pi_{\sF_\Gamma}$ for all  $\la$  in $\BD$. Here
$\Pi_{\sF_\Gamma}$ is the orthogonal projection from $\sD_\Gamma$
onto $\sF_\Gamma$. We claim that for the given contraction
$\Gamma$, the set of all Schur pairs $\{F,G\}$ associated with the
data set $\{A,T^\prime,U^\prime,R,Q\}$ and satisfying
(\ref{fundeq}) is parameterized by the set
$\eS_\Omega(\sD_\Gamma,\sD_\Gamma)$. To make this precise, we
first define a mapping $J_\Gamma$ from
$\eS(\sD_\Gamma,\sD_\Gamma)$ into $\eS(\sD_A, \sD_{T'}\oplus
\sD_A) $ as follows
\begin{equation}\label{defJ}
J_\Gamma C = \left[ \begin{array}{c}
  F\\
  G
\end{array}\right], \quad \left[ \begin{array}{c}
  F(\lambda) \\
   \noalign{\vskip4pt}
  G(\lambda)
\end{array}\right] =
\left[ \begin{array}{c}
  2 \Theta(\lambda) \left(W(\lambda) + I\right)^{-1}\\
  \noalign{\vskip4pt}
  \lambda^{-1} \left(W(\lambda) - I\right)\left(W(\lambda) + I\right)^{-1}
\end{array}\right],
\end{equation}
where
\begin{eqnarray}
\Theta(\lambda)d&=&(\Gamma d)(\lambda), \quad d\in \sD_A,\nonumber\\
 \noalign{\vskip4pt}
W(\lambda) &=& \Gamma^*(I+\la S^*)(I-\la S^*)^{-1}\Gamma
+D_\Gamma(I+\la C(\la))(I-\la C(\la))^{-1}D_\Gamma,\quad\la\in\BD.
\label{defW}
\end{eqnarray}
Here $S$ is the unilateral shift on $H^2(\sD_{T^\prime})$ and $E$
is the canonical embedding of $\sD_{T^\prime}$ onto the set of
constant function in $H^2(\sD_{T^\prime})$. We are now ready to
state the second main theorem.

\begin{theorem}\label{th2} Let $B$ in $(\ref{defB})$ be
a solution to the RCL problem for the data set
$\{A,T^\prime,U^\prime,R,Q\}$ with $U^\prime$ being given by
$(\ref{szns})$, and let $\Gamma$ be the contraction determined by
$B$ via $(\ref{defB})$. Then the mapping $J_\Gamma$ from
$\eS(\sD_\Gamma,\sD_\Gamma)$ into $\eS(\sD_A,\sD_{T^\prime} \oplus
\sD_A )$ defined in $(\ref{defJ})$  maps
$\eS_\Omega(\sD_\Gamma,\sD_\Gamma)$ in a  one to one way onto the
set of all Schur pairs $\{F,G\}$ associated with the given data
set such that $(\ref{sols1})$ and $(\ref{sols2})$ hold.
\end{theorem}

To give some further insight in the set
$\eS_\Omega(\sD_\Gamma,\sD_\Gamma)$ appearing in Theorem
\ref{th2}, put $\sG_\Gamma=\sD_\Gamma\ominus\sF_\Gamma$, and let
$\Pi_{\sF_\Gamma}$ and $\Pi_{\sG_\Gamma}$ be the orthogonal
projections from $\sD_\Gamma$ onto $\sF_\Gamma$ and $\sG_\Gamma$,
respectively. Using Corollary XXVII.5.3 in \cite{ggk2} it follows
that $C\in\eS_\Omega(\sD_\Gamma,\sD_\Gamma)$ if and only if
\begin{equation}\label{CC1}
C(\la)=\Omega\Pi_{\sF_\Gamma}+D_{\Omega^*}C_1(\la)\Pi_{\sG_\Gamma},\quad
\la\in\BD,
\end{equation}
for some function $C_1$ in the Schur class
$\eS(\sG_\Gamma,\sD_{\Omega^*})$. Moreover, $C$ and $C_1$ in
$(\ref{CC1})$ determine each other uniquely. Hence, instead of
$\eS_\Omega(\sD_\Gamma,\sD_\Gamma)$, we can say, in Theorem
$\ref{th2}$, that the set of all Schur pairs $\{F,G\}$ satisfying
$(\ref{sols2})$ correspond to $\eS(\sG_\Gamma,\sD_{\Omega^*})$ in
a one to one way.

A similar remark applies to the set of Schur pairs appearing in
Theorem \ref{thrc}. To see this, notice that a pair of functions
$\{F,G\}$ is a Schur pair associated to the data set
$\{A,T^\prime,U^\prime,R,Q\}$ if and only if
\[
\col[F,G]\in \{H\in\eS(\sD_A,\sD_{T'}\oplus\sD_A)\, :\,
H(\la)|\sF=\om\mbox{ for all }\la\in\BD\}.
\]
Therefore, the set of Schur pairs associated to the given data set
is in one to one correspondence to $\eS(\sG,\sD_{\om^*})$, where
$\sG=\sD_A\ominus\sF$.

\medskip We conclude this section with
the commutant lifting theorem as given by Theorem XIII.3.4 in
\cite{ff}. We show how this result can be derived from Theorems
\ref{thrc} and \ref{th2}.

\begin{theorem}\label{CLP}
Let $\{A,T^\prime,U^\prime,R,Q\}$ be a data set with $U^\prime$
being given by $(\ref{szns})$, $\sH_0=\sH$, $R=I_\sH$ and $Q$ an isometry on
$\sH$. Then all solutions to the
corresponding RCL problem are given by
\begin{equation}\label{classols1.1}
B=\mat{c}{A\\\Gamma D_A}:\sH\to\mat{c}{\sH'\\H^2(\sD_{T'})},
\end{equation}
where $\Gamma$ is a contraction from $\sD_A$ into $H^2(\sD_{T'})$ given by
\begin{equation}\label{classols1.2}
(\Gamma d)(\lambda)=F(\lambda)(I-\lambda G(\lambda))^{-1}d,\quad
d\in\sD_A,\lambda\in\BD,
\end{equation}
with $\{F,G\}$ an arbitrary Schur pair associated with the given
data set. The solution $B$ and the Schur pair $\{F,G\}$ in
$(\ref{classols1.1})$ and $(\ref{classols1.2})$ determine each
other uniquely. Finally, there exists only one solution to the
given RCL problem if and only if $\sF=\sD_A$ or
$\om\sF=\sD_{T'}\oplus\sD_A$.
\end{theorem}

\bpr The representation of all solutions follows immediately from
Theorem \ref{thrc}. Obviously, the Schur pair $\{F,G\}$ in
(\ref{classols1.2}) determines $B$ uniquely. To prove the converse
implication, let $B$ be a  solution to the corresponding RCL
problem for the given data set, and let $\Gamma$ be the
contraction from $\sD_A$ into $H^2(\sD_{T'})$ given by
(\ref{classols1.1}). By Theorem \ref{th2} it suffices to show that
the set $\sS_\Omega(\sD_\Gamma,\sD_\Gamma)$ consists of one
element only. Recall that in the commutant lifting setting
$R^*R=Q^*Q$, and hence, as has been remarked in the paragraph
preceding (\ref{sxi}), in this case the operator $\Omega$ is an
isometry. Moreover, from the definition of $\Omega$ we obtain that
\[
\im\Omega=\overline{D_\Gamma D_AR\sH_0}=\overline{D_\Gamma D_A
\sH} =\sD_\Gamma.
\]
Thus $\Omega$ is a unitary operator from $\sF_\Gamma$ onto
$\sD_\Gamma$, and hence  $\sD_{\Omega^*}=\{0\}$. But then the
remark made in the first paragraph after Theorem \ref{th2} shows
that $\sS_\Omega(\sD_\Gamma,\sD_\Gamma)$ is a singleton.

Now that we know that every solution uniquely corresponds to a
Schur pair, we see that there is only one solution if and only if
there is only one corresponding Schur pair. {}From the remark made
in the second paragraph after Theorem \ref{th2} we see that  the
latter happens if and only if $\eS(\sG,\sD_{\om^*})$ consists of
the zero element only. In other words, there exists a unique
solution if and only if $\sG=\{0\}$ or $\sD_{\om^*}=\{0\}$. From
$\sG=\sD_A\ominus\sF$ it follows that $\sG=\{0\}$ if and only if
$\sF=\sD_A$. Since in the commutant lifting setting the operator
$\omega$ is an isometry (see the paragraph containing
(\ref{altfundeq})), we have
$\im\om=\ker\om^*=(\sD_{T'}\oplus\sD_A)\ominus\sD_{\om^*}$. Hence
$\sD_{\om^*}=\{0\}$ is equivalent to
$\om\sF=\sD_{T'}\oplus\sD_A$.\epr

\medskip
For the commutant lifting setting representations of all solutions
by formulas of the type (\ref{classols1.2}) date back to
\cite{ACF80}, see also \cite{FF2}. The proofs of Theorems \ref{thrc} and
\ref{th2} will be given in the third section.

\section{Operator-valued positive real functions and harmonic majorants}
\setcounter{equation}{0}

Let $\Theta$ be a $\eL(\sE,\sY)$-valued analytic function on
$\mathbb{D}$, where $\sE$ and $\sY$ are Hilbert spaces.
 We say that $\Theta$ belongs to $H^2(\eL(\sE,\sY))$ if for
each $a \in \sE$ the function $\Theta(\cdot)a$ belongs to
$H^2(\sY)$. The latter condition is equivalent to the requirement
that $\sum_{\nu=0}^\infty \|\Theta_\nu a\|^2 < \infty$ for all $a$
in $\sE$. Here and in the sequel $\Theta_0,\Theta_1,\Theta_2,
\dots$ are the Taylor coefficients of $\Theta$ at zero. If
$\Theta$ is in $H^2(\eL(\sE,\sY))$, then $\Theta$ uniquely defines
an operator $\Gamma$ from $\sE$ into $ H^2(\sY)$ by
\begin{equation}  \label{defGa}
(\Gamma a)(\lambda)  = \Theta(\lambda) a, \quad a\in\sE,
\lambda \in \mathbb{D}.
\end{equation}
In this case, we say that $\Gamma$ is the {\em operator
associated} with $\Theta$. On the other hand, if $\Gamma$ is an
operator mapping $\sE$ into $H^2(\sY)$, then the relation
$\Theta(\lambda) a = (\Gamma a)(\lambda)$ for $a$ in $\sE$ and
$\lambda$  in $\mathbb{D}$ uniquely defines a function $\Theta$ in
$H^2(\eL(\sE,\sY))$. In this case, we say with a slight abuse of
terminology that $\Theta$ is the \emph{symbol}  of $\Gamma$.

As before, let $\Theta$ be a function in $H^2(\eL(\sE,\sY))$, and
let $\Gamma$ be the operator associated with $\Theta$. Throughout this
section $S$ is the block forward shift on $H^2(\sY)$, and $E$ the
canonical embedding from $\sY$ onto the constant functions in
$H^2(\sY)$, that is, $(E y)(\lambda) \equiv y $ on $\opendisc$. In
this  case, $\Theta_n = E^* (S^*)^n\Gamma$ for all non-negative
integers $n $. Hence $\Theta$ admits a state space realization  of
the following form:
\begin{equation}\label{reprK}
\Theta(\lambda) = E^* (I-\lambda S^*)^{-1}\Gamma, \quad \lambda
\in\mathbb{D}.
\end{equation}

With $\Theta$ as above we associate the $\eL(\sE,\sE)$-valued function
\begin{equation}\label{defV}
V(\lambda) =
\Gamma^*\Gamma + 2\lambda\Gamma^* (I-\lambda S^*)^{-1}S^*\Gamma,
\quad  \lambda \in\mathbb{D},
\end{equation}
where $\Gamma$ is the operator associated with $\Theta$ via (\ref{reprK}).
An easy computation shows that $V$ can also be written as
\begin{equation}\label{defVb}
V(\lambda)=\Gamma^*(I+\lambda S^*)(I-\lambda S^*)^{-1}\Gamma,\quad
\lambda\in\BD.
\end{equation}
Obviously, $V$ is analytic on $\mathbb{D}$. Using $EE^*=I-SS^*$,
we see from (\ref{reprK}) and (\ref{defV}) that the Taylor
coefficients $\{V_n\}_0^\infty$ of $V$ at zero are given by
\[
V_0 = \Gamma^*\Gamma = \sum_{\nu=0}^\infty \Theta_\nu^*\Theta_{\nu}
\quad \mbox{and} \quad
 V_n = 2\Gamma^* S^{*n} \Gamma = 2\sum_{\nu=0}^\infty
\Theta_\nu^*\Theta_{\nu +n}, \quad\mbox{for all } n \geq 1.
\]
The results below show that $V$ is positive real,
and therefore we shall refer to $V$ as the \emph{positive real
function defined by} $\Theta$.

Recall that a $\eL(\sE,\sE)$-valued function $W$ is \emph{positive
real} if $W$ is analytic on $\opendisc$ and
\[
\Re W(\lambda)=\half\big(W(\lambda)^*+W(\lambda)\big)\geq 0,
\quad \lambda \in\mathbb{D}.
\]
It is known (see, e.g., \cite{fk02a}, Section 1.2) that a
$\eL(\sE,\sE)$-valued function $W$  which is analytic at zero,
$W(\lambda)=\sum_{\nu=0}^\infty \lambda^\nu W_\nu$ say, is positive
real if and only if for each $n$ the $n\times n$ Toeplitz operator matrix
$T_{\Re W,\, n}$ given  by
\begin{equation}\label{mattoep}
T_{\Re W,\, n}=
\half \left[
\begin{array}{cccc}
 W_0^*+W_0 & W_1^* & \cdots & W_{n-1}^*\\
W_1& W_0^*+W_0  &\cdots  &W_{n-2}^* \\
  \vdots &\vdots  & \ddots & \vdots \\
W_{n-1}&W_{n-2}&\cdots   & W_0^*+W_0 \\
\end{array}\right],
\end{equation}
defines a non-negative operator on $\sE^n$.

Our aim in this section is to prove the following theorem which
can be viewed as an operator valued version of a classical result
on harmonic majorants, cf., Section 2.6 in \cite{D70}.

\begin{theorem}\label{param}
Let $\Theta$ be a function in $H^2(\eL(\sE,\sY))$ such that the associated
operator $\Gamma$ is a contraction from $\sE$ into $H^2(\sY)$.
The set of all positive real functions $W$ with values in
$\eL(\sE,\sE)$ satisfying
\begin{equation}\label{eq:bound1}
\Theta(\lambda)^* \Theta(\lambda) \leq\Re W(\lambda) \quad
\mbox{for all } \lambda\in\mathbb{D} \mbox{ and } W(0) = I
\end{equation}
is parameterized by $\eS(\sD_\Gamma,\sD_\Gamma)$. More
precisely, all positive real functions $W$ on $\mathbb{D}$
satisfying  $(\ref{eq:bound1})$  are given by
\begin{equation}\label{eq:para}
W(\lambda) = V(\lambda)+ D_\Gamma \left(I + \lambda
C(\lambda)\right) \left(I - \lambda C(\lambda)\right)^{-1}
D_\Gamma, \quad \lambda\in\mathbb{D},
\end{equation}
where $V$ on $\mathbb{D}$ is given by $(\ref{defV})$, and $C$ is
an arbitrary function in $\eS(\sD_\Gamma,\sD_\Gamma)$. Moreover,
$W$ and $C$ in $(\ref{eq:para})$ determine each other uniquely.
Finally, there is only one positive real function $W$ satisfying
$(\ref{eq:bound1})$ if and only if $\Gamma$ is an isometry. In
this case $W = V$ is the only function satisfying
$(\ref{eq:bound1})$.
\end{theorem}

In order to prove the above theorem it will be convenient to first
prove a lemma and to review some theory concerning the Cayley
transform of operator-valued functions.

\begin{lemma} \label{lemineq}
Let $\Theta\in H^2(\eL(\sE,\sY))$, and $V$  be the
$\eL(\sE,\sE)$-valued function defined by $(\ref{defV})$. Then $V$
is positive real. More precisely,
\begin{equation}\label{basicineq}
\Theta(\lambda)^*\Theta(\lambda)\leq \Re V(\lambda),
\quad  \lambda \in \mathbb{D}.
\end{equation}
Furthermore, if $W$ is any $\eL(\sE,\sE)$-valued positive real
function  such that $\Theta(\lambda )^*\Theta(\lambda ) \leq \Re
W(\lambda)$ for all $\lambda \in\mathbb{D}$, then $W-V$ is
positive real.
\end{lemma}

To give some further insight in (\ref{basicineq}), let us consider the
scalar case, that is, $\sE$ and $\sY$ are equal to $\BC$. In that case
formula (\ref{defVb}) can be rewritten as
\[
V(\la)=\frac{1}{2\pi}\int_0^{2\pi}\frac{e^{i\om}+\la}{e^{i\om}-\la}
|\theta(e^{i\om})|^2 d\om,\quad\la\in\BD,
\]
and the above lemma is well known (see the proof of Theorem 2.12 in
\cite{D70}).
In fact, in the scalar case $\Re V$ is known as the least harmonic
majorant of $|\theta(\cdot)|^2$.

\smallskip
\noindent {\bf Proof of Lemma~\ref{lemineq}.}
We split the proof into three parts. In the
first part we prove (\ref{basicineq}).

\smallskip\noindent\textit{Part 1.} Take $\la\in\BD$.
For convenience set
$\Phi(\lambda) = (I - \lambda S^*)^{-1}$. Using (\ref{reprK}) and
(\ref{defV}), we have
\[
\Theta(\lambda) = E \Phi(\lambda)\Gamma
\quad \mbox{and} \quad
V(\lambda)=\Gamma^*\Gamma
+2\lambda\Gamma^*\Phi(\lambda)S^*\Gamma.
\]
Note that $\Phi(\lambda)=I+\lambda\Phi(\lambda)S^*$. Since
$E^* E  + SS^* = I$, we obtain
\begin{eqnarray*}
\Theta(\lambda)^*\Theta(\lambda) &=& \Gamma^*\Phi(\lambda)^* E^* E \Phi(\lambda)\Gamma =
\Gamma^*\Phi(\lambda)^*(I - SS^*)\Phi(\lambda )\Gamma \\\noalign{\vskip4pt} &=&
\Gamma^*\Phi(\lambda )^*\Phi(\lambda )\Gamma -
  \Gamma^*\Phi(\lambda )^*SS^*\Phi(\lambda )\Gamma \\
\noalign{\vskip4pt} &=& \Gamma^* [I+\bar{\lambda}S \Phi(\lambda)^* ] [I+ \lambda
\Phi(\lambda) S^* ]\Gamma -
 \Gamma^*\Phi(\lambda )^*SS^*\Phi(\lambda )\Gamma\\
\noalign{\vskip4pt} &=& \Gamma^*\Gamma+ \bar{\lambda } \Gamma^*S\Phi(\lambda)^*
\Gamma +\lambda
\Gamma^*\Phi(\lambda)S^*\Gamma +\\
\noalign{\vskip4pt} &&\hspace{2.5cm}
+|\lambda|^2\Gamma^*\Phi(\lambda)^*SS^*\Phi(\lambda)\Gamma
-\Gamma^*\Phi(\lambda )^*SS^*\Phi(\lambda )\Gamma\\
\noalign{\vskip4pt} &=&  \half\left(V(\lambda) +
V(\lambda)^*\right)- (1-|\lambda|^2)\Gamma^*(I - \bar{\lambda} S)^{-1}SS^*(I -
\lambda  S^*)^{-1}\Gamma.
\end{eqnarray*}
The last term is non-negative. Thus (\ref{basicineq}) holds. In
particular, $V$ is positive real.

\smallskip\noindent\textit{Part 2.}
Fix $0<r<1$, and set $\widetilde{\Theta}(z) = \Theta(r z)$
for each $z \in \opendisc$. Notice that $\widetilde{\Theta}$ is analytic in open
neighborhood of $\closedisc$, the closure of the open unit disc
$\opendisc$. Let $\widetilde{\Gamma}$ be the operator from $\sE$ into
$H^2(\sY)$  associated with
$\widetilde{\Theta}$, that is,
$(\widetilde{\Gamma} a)(z) = \widetilde{\Theta}(z)a$ for
$a\in \sE$ and $z \in \mathbb{D}$. Thus $\widetilde{\Gamma} = \Lambda_r\Gamma$,
where $\Lambda_r$ is the   operator on $H^2(\sY)$ defined by
\[
(\Lambda_r h)(z)  = h(r z), \quad h \in H^2(\sY) , z
\in \mathbb{D}.
\]
Note that $\Lambda_r$ is bounded and $\lim_{r\uparrow1}\Lambda_r=I$ with
pointwise convergence.
Let $\widetilde{V}$ be the positive real function defined by
$\widetilde{\Theta}$. Thus
\[
\widetilde{V}(\lambda) = \Gamma^*\Lambda_r^2\Gamma+2\lambda\Gamma^* \Lambda_r
(I-\lambda S^*)^{-1}S^*\Lambda_r\Gamma, \quad \lambda \in\mathbb{D}.
\]
Since $\Lambda_rS=rS\Lambda_r$, we have $\Lambda_r(I-\lambda S)^{-1}=(I-\lambda r
S)^{-1}\Lambda_r$ for each $\lambda\in \opendisc$.
Taking adjoints and replacing $\la$ by $\bar{\la}$ we see that
$(I-\la S^*)^{-1}S^*\Lambda_r=r\Lambda_r(I-\la rS^*)^{-1}S^*$ and hence
$\widetilde{V}$ is also analytic on an open neighborhood of $\closedisc$.

{}From the first part of the proof we know that for each $\lambda$ in
$\opendisc$ we have
\begin{eqnarray*}
\Re \widetilde{V}(\lambda)-\widetilde{\Theta}(\lambda )^*\widetilde{\Theta}(\lambda
)&=&(I-|\lambda|^2)\Gamma^*\Lambda_r(I - \bar{\lambda} S)^{-1}SS^*(I - \lambda
S^*)^{-1}\Lambda_r\Gamma\\
\noalign{\vskip4pt} &=& (I-|\lambda|^2)\Gamma^*(I - \bar{\lambda}r
S)^{-1}\Lambda_r S S^*\Lambda_r(I - \lambda rS^*)^{-1}\Gamma.
\end{eqnarray*}
Since all functions involved are analytic on an open neighborhood
of $\closedisc$, we conclude that
\[
\widetilde{\Theta}(e^{\imath \omega}
)^*\widetilde{\Theta}(e^{\imath \omega } )= \Re
\widetilde{V}(e^{\imath \omega }), \quad 0 \leq \omega\leq 2\pi.
\]

Let $W$ be a  positive real function with values in
$\eL(\sE,\sE)$ such that $\Theta(\lambda )^*\Theta(\lambda ) \leq
\Re W(\lambda)$ for all $\lambda \in\mathbb{D}$. Set
$\widetilde{W}(\lambda)=W(r\lambda)$ for each $\lambda
\in\mathbb{D}$. Then $\widetilde{\Theta}(\lambda
)^*\widetilde{\Theta}(\lambda ) \leq \Re \widetilde{W}(\lambda)$
for all $\lambda \in\mathbb{D}$. Again $\widetilde{W}$ is analytic
on an open neighborhood of $\closedisc$, and thus, by continuity,
$\widetilde{\Theta}(e^{\imath \omega }
)^*\widetilde{\Theta}(e^{\imath \omega } ) \leq \Re
\widetilde{W}(e^{\imath \omega })$ for each $0\leq \omega\leq
2\pi$. But then we can use the result of the previous paragraph to
show that
\begin{equation}\label{ineq0}
\Re \widetilde{V}(e^{\imath \omega })\leq \Re
\widetilde{W}(e^{\imath \omega}),  \quad 0\leq \omega\leq 2\pi.
\end{equation}

Next we show that the latter inequality implies that
$\widetilde{W}-\widetilde{V}$ is positive real. To accomplish this, let
$L_{\Re \widetilde{V}}$ and $L_{\Re \widetilde{W}}$ be the block
Laurent operators on $\ell^2(\sY)$ defined by $\Re \widetilde{V}$
and $\Re \widetilde{W}$, respectively. Since $\Re \widetilde{V}$
and $\Re \widetilde{W}$ are both continuous on the unit circle
$\mathbb{T}$, these operators are well defined and bounded.
Furthermore, the inequality (\ref{ineq0}) implies that $L_{\Re
\widetilde{V}}\leq L_{\Re \widetilde{W}}$. Taking the compression
to $\ell_+^2(\sY)$ this implies that  $T_{\Re \widetilde{V}}\leq
T_{\Re \widetilde{W}}$, where $T_{\Re \widetilde{V}}$ and $T_{\Re
\widetilde{W}}$ are the block Toeplitz operators on
$\ell_+^2(\sY)$ defined by $\Re \widetilde{V}$ and $\Re
\widetilde{W}$, respectively. Next, taking an $n$-th section of
these block Toeplitz operators, we obtain that $T_{\Re
\widetilde{V},\,n}\leq T_{\Re \widetilde{W},\,n}$ for all integers
$n \geq 0$. This implies (see the paragraph before Lemma
\ref{lemineq}) that $\widetilde{W}-\widetilde{V}$ is positive
real.

\smallskip\noindent\textit{Part 3.}
We continue to use the notation introduced in the preceding part,
but now we make the dependence on the parameter $r$ explicit. Thus
for $\widetilde{V}$  we write  $V_{(r)}$, and for $\widetilde{W}$
we write  $W_{(r)}$. Define
\[
\Delta = W - V, \quad \Delta_{(r)}= W_{(r)} - V_{(r)},\quad
\mbox{for each }0<r<1.
\]
The result of the previous part shows that $\Delta_{(r)}$ is
positive real for each $0 < r < 1$. Furthermore, for $r\uparrow 1$
the $n$-th Taylor coefficient of $\Delta_{(r)}$ converges
pointwise (i.e., in the strong operator topology) to the
$n$-th Taylor coefficient of $\Delta$. Here $n$ is an arbitrary
non-negative integer. Hence for each $n=0,1,2,\dots$ we see that
$T_{\Re \Delta_{(r)}, \, n}x$ converges to $T_{\Re \Delta, \,
n}x$ for each $x\in\sE^n$ as $r\uparrow 1$. Since the operators
$T_{\Re \Delta_{(r)}, \, n}$ are non-negative, the same holds true
for $T_{\Re \Delta, \, n}$. This shows that $\Delta = W - V$ is
positive real. \epr

\medskip \noindent\textbf{Positive real functions and the Cayley
transform.}  For $C$ in $\eS(\sE,\sE)$ consider the map
\begin{equation}\label{Ctransf}
C\mapsto \K, \quad \mbox{where} \quad \K(\lambda) = \left(I +
\lambda C(\lambda)\right) \left( I-\lambda C(\lambda)\right)^{-1}
\quad \mbox{for all }\lambda \in \mathbb{D}.
\end{equation}
Since $C(\lambda)$ is contractive for each $\lambda\in
\mathbb{D}$, the function $K$ is well defined by (\ref{Ctransf}).
The map $C\mapsto K$ in (\ref{Ctransf}) establishes a one to  one
correspondence between the Schur class $\eS(\sE,\sE)$ and the set
of all positive real functions $\K$ satisfying $\K(0) = I$.
Indeed, if $\K$ is defined by (\ref{Ctransf}) for some $C\in
\eS(\sE,\sE)$, then $\K$ is analytic in $\mathbb{D}$ and  $\K(0) =
I$ while
\begin{equation}\label{eq:realpart}
\Re \K(\lambda) = \left(I - \lambda C(\lambda)\right)^{-*} \left(I
- |\lambda|^2 C(\lambda)^*C(\lambda)\right) \left(I - \lambda
C(\lambda)\right)^{-1}, \quad \lambda\in\mathbb{D}.
\end{equation}
It follows that $\Re \K(\lambda)>0$ for each $\lambda \in
\mathbb{D}$, and hence $\K$ is positive real. Conversely, for a
positive real function $\K$ satisfying $\K(0) = I$, the function
$C$ given by
\begin{equation}\label{eq:invCayley}
C(\lambda) = \frac{1}{\lambda}\left(\K(\lambda) - I\right)\left(I
+ \K(\lambda)\right)^{-1}, \quad 0\not = \lambda\in\mathbb{D},
\end{equation}
is well defined and belongs to $\eS(\sE,\sE)$.

If $C$ belongs to $\eS(\sE,\sE)$, then we call $K$ defined by
(\ref{Ctransf}) the \emph{Cayley transform} of $C$. If $K$ is
positive real with $K(0) = I$, then $C$ defined by
(\ref{eq:invCayley}) will be called the \emph{inverse Cayley
transform} of $K$.

\smallskip
\noindent {\bf Proof of Theorem~\ref{param}.}  Let $C$ be a function in
$\eS(\sD_\Gamma,\sD_\Gamma)$,
 and define $W$ by (\ref{eq:para}).
Then $W (\lambda) = V(\lambda) + D_\Gamma \K (\lambda)D_\Gamma$
for each $\lambda$  in $\mathbb{D}$, where $\K$ on $\mathbb{D}$ is
the Cayley transform of $C$. Note that $V(0) = \Gamma^*\Gamma$.
Hence $W(0) = V(0)+ I -\Gamma^* \Gamma = I$. By consulting  Lemma
\ref{lemineq}, we have
\[
\Re W(\lambda)=\Re V(\lambda) + D_\Gamma (\Re \K(\lambda))D_\Gamma
\geq \Re V(\lambda) \geq \Theta(\lambda)^*\Theta(\lambda)  \geq 0,
\quad \lambda\in\mathbb{D}.
\]
Therefore  $W$ is a positive real function satisfying
(\ref{eq:bound1}).

Conversely, assume that $W$ is a positive real function satisfying
(\ref{eq:bound1}). According to Lemma \ref{lemineq},  we have that
$\Re W(\lambda) \geq \Re V(\lambda)$ for all $\lambda$ in
$\mathbb{D}$. Hence, the function $\Delta = W - V$ is a positive
real function on $\mathbb{D}$ that satisfies $\Delta(0) = W(0) -
V(0) = I - \Gamma^*\Gamma=D_\Gamma^2$.

We claim that  $\Delta$ admits a unique factorization of the form
$\Delta(\lambda) = D_\Gamma \K(\lambda) D_\Gamma$, where $\K$ is a
positive real function  with values in
$\eL(\sD_\Gamma,\sD_\Gamma)$ and $\K(0) = I$. To see this let
$\{\Delta_n\}_0^\infty$ be the Taylor coefficients of $\Delta$ at
the origin. Since $T_{\Re \Delta,n}$ is a positive Toeplitz matrix
and $\Delta(0) = D_\Gamma^2$, we see that
\begin{equation}\label{pos}
 \left[ \begin{array}{cc}
  2 D_\Gamma^2 & \Delta_n^* \\
  \Delta_n & 2 D_\Gamma^2 \\
\end{array}\right] \geq 0, \quad n=0,1,2,\ldots.
\end{equation}
Recall (see Theorem XVI.1.1. in \cite{ff}) that a $2\times 2$
operator matrix
\[
 \left[ \begin{array}{cc}
  A & B^* \\
  B &  A \\
\end{array}\right] \mbox{ on } \left[ \begin{array}{cc}
  \sE \\
  \sE \\
\end{array}\right]
\]
induces a positive operator on $\sE\oplus\sE$ if and only if $B =
A^{1/2} \Phi A^{1/2}$ for some contraction $\Phi$ on $\overline{A
\sE}$. In this case, $B$ and $\Phi$ uniquely determine each other.
So from (\ref{pos}) we see that there exists a unique operator
$\K_n$ on $\sD_\Gamma$ such that $\Delta_n = D_\Gamma \K_n
D_\Gamma$ for all integers $n \geq 0$, and $\K_0 = I$. Let $T_{\Re
\K,\,n}$ the $n\times n$ block Toeplitz operator matrix obtained
by replacing $W_j$ by $\K_j$ in (\ref{mattoep}).
Notice  $\eD_n^*  T_{\Re \K,\,n} \eD_n = T_{\Re \Delta,\,n}$,
where $\eD_n$ is the diagonal operator matrix
$\mbox{diag} \{D_\Gamma\}_1^n$ acting on
$\oplus_1^n \sD_\Gamma$.
Since $T_{\Re \Delta,\, n}$ is positive, and $\eD_n$
is onto a dense set in $\oplus_1^n \sD_\Gamma$, it follows that
$T_{\Re K,\, n}$ is positive for each integer $n \geq 0$.
Hence  $\K(\lambda) = \sum_{n=0}^\infty \lambda^n \K_n$ is a
positive real function. Therefore
$\Delta(\lambda) = D_\Gamma \K(\lambda) D_\Gamma$
where $\K$ is a positive real function  satisfying $\K(0) = I$,
which proves our claim.

Let $C$ on $\mathbb{D}$ be the inverse Cayley transform of $\K$.
Then $C$ is a function in
$\eS(\sD_\Gamma,\sD_\Gamma)$, and we have
\[
\K(\lambda) = \left(I + \lambda C(\lambda)\right) \left(I -
\lambda C(\lambda)\right)^{-1}, \quad \lambda\in\mathbb{D}.
\]
Hence $W$ is given by (\ref{eq:para}) with $C\in
\eS(\sD_\Gamma,\sD_\Gamma)$ being the inverse Cayley
transform of the positive real function $\K$ uniquely determined
by $\Delta(\lambda) = D_\Gamma \K(\lambda) D_\Gamma$. Recall  that
the inverse Cayley transform is a bijective mapping from the set
of positive real functions $\K$ with $\K(0) = I$ onto
$\eS(\sD_\Gamma,\sD_\Gamma)$. Thus   $\K$ and $C$
uniquely determine each other.

Moreover, since $\Delta$ and $K$ determine each other uniquely and
the Cayley transform is bijective, we obtain that $C$ and
$W$ in (\ref{eq:para}) determine each other uniquely. \epr

\section{Proofs of the main theorems}
\setcounter{equation}{0}

In this section we proof Theorems \ref{thrc} and \ref{th2}.
Throughout this section $\{A,T^\prime,U^\prime,R,Q\}$ is a fixed
data set with $U'$ being given by~(\ref{szns}). As mentioned in
Section~\ref{sec:mainth}, an operator $B$ from $\sH$ into
$\sH'\oplus H^2(\sD_{T'})$ is a solution to the corresponding RCL
problem if and only if $B$ admits a representation of the form
\begin{equation}\label{eq:Bdec}
B=\mat{c}{A\\\Gamma D_A}:\sH\to\mat{c}{\sH'\\H^2(\sD_{T'})},
\end{equation}
with $\Gamma$ a contraction from $\sD_A$ into $H^2(\sD_{T'})$ satisfying
\begin{equation}\label{altfundeq2}
E\om_1+S\Gamma\om_2=\Gamma|\sF.
\end{equation}
Here $S$ denotes the unilateral shift on $H^2(\sD_{T'})$ and $E$ is
the canonical embedding of $\sD_{T'}$ onto the space of constant functions in
$H^2(\sD_{T'})$ defined by $(E d)(\lambda)\equiv d$ for all $d\in \sD_{T'}$.

As a first step towards the proofs of Theorem \ref{thrc} and \ref{th2} it will
be convenient first to consider the case when the space $\sF$ in
(\ref{altfundeq2}) consists of the zero element only. In that case the only
constraint on the operator $\Gamma$ in (\ref{eq:Bdec}) is that it has to be a
contraction. It follows that for $\sF=\{0\}$ our two main theorems reduce to
the following result.

\begin{theorem}\label{th:H2}
Let $\Gamma$ be an operator from $\sE$ into $H^2(\sY)$. Then $\Gamma$ is a
contraction if and only if $\Gamma$ admits a representation of the form
\begin{equation}\label{Garep}
(\Gamma e)(\lambda)=F(\lambda)(I-\lambda G(\lambda))^{-1}e,\quad
e\in\sE, \lambda\in\BD,
\end{equation}
where $\col[F,G]$ is any function in $\eS(\sE,\sY\oplus\sE)$.
Moreover, if $\Gamma$ is a contraction, then there is a one to one
correspondence between $\eS(\sD_\Gamma,\sD_\Gamma)$ and the set of
all Schur class functions $\col[F,G]$ in $\eS(\sE,\sY\oplus\sE)$,
that satisfy $(\ref{Garep})$. To be precise, let $J_\Gamma$ be the
map from $\eS(\sD_\Gamma,\sD_\Gamma)$ into $\eS(\sE,\sY\oplus\sE)$
defined by
\begin{equation}\label{JGaabs}
J_\Gamma C=\mat{c}{F\\G}\ (C\in\eS(\sD_\Gamma,\sD_\Gamma)),\
\mbox{where}\
\mat{c}{F(\lambda)\\G(\lambda)}=\mat{c}{2\Theta(\lambda)(W(\la)+I)^{-1}\\
\lambda^{-1}(W(\lambda)-I)(W(\la)+I)^{-1}}
\end{equation}
with $\Theta$ the symbol of $\Gamma$, see $(\ref{defGa})$,
and
\begin{equation}\label{defW2}
W(\la)=\Gamma^*(I+\la S^*)(I-\la S^*)^{-1}\Gamma +D_\Gamma(I+\la
C(\la))(I-\la C(\la))^{-1}D_\Gamma,\quad \la\in\BD.
\end{equation}
Then $J_\Gamma$ is a one to one mapping from
$\eS(\sD_\Gamma,\sD_\Gamma)$ onto the set of all functions
$\col[F,G]$ in $\eS(\sE,\sY\oplus\sE)$  that satisfy
$(\ref{Garep})$. In particular, the representation in
$(\ref{Garep})$ is unique if and only if $\Gamma$ is an isometry.
\end{theorem}

In a somewhat different, less explicit form, Theorem \ref{th:H2}
appears in the introduction of \cite{FtHK04}, see Corollaries 0.3
and 0.4 in \cite{FtHK04}. These corollaries were obtained as
immediate consequences of the description of all solutions to the
relaxed commutant lifting problem given in \cite{FtHK04}. In the
present paper we follow a different direction: we first proof
Theorem \ref{th:H2}, and then derive Theorems \ref{thrc} and
\ref{th2} as further refinements of Theorem \ref{th:H2}.

Theorem \ref{th:H2} has other partial predecessors in the
literature. For example, when $\sE=\sY=\BC$ and $\Gamma$ is an
isometry, the representation (\ref{Garep}) immediately follows
from the description of $H^2$ functions of unit norm given in
\cite{Sar89}, page 490. When $\sE=\BC^q$ and $\sY=\BC^p$ the first
statement in Theorem \ref{th:H2} is Theorem 2.2 in \cite{ABP95}.
The second and third part of Theorem \ref{th:H2} seem to be new,
even in the scalar case.

\smallskip\noindent
{\bf Proof of Theorem \ref{th:H2}.} Let $\Theta$ be the symbol of
$\Gamma$, see (\ref{defGa}). Take for $C$ any function in
$\eS(\sD_\Gamma,\sD_\Gamma)$, and define functions $F$ and $G$ by
(\ref{JGaabs}) and (\ref{defW2}). Then $F$ is a
$\eL(\sE,\sY)$-valued function and $G$ is a $\eL(\sE,\sE)$-valued
function. {}From Theorem \ref{param} we obtain that $W$ in
(\ref{defW2}) is a positive real function satisfying
(\ref{eq:bound1}). Note that $G$ is the inverse Cayley transform
of $W$. Hence $G$ is a function in $\eS(\sE,\sE)$. Moreover, for
each $\la\in\BD$ we have
\begin{eqnarray}
I-\la G(\la)
&=&I-(W(\la)-I)(W(\la)+I)^{-1}
=\left((W(\la)+I)-(W(\la)-I)\right)(W(\la)+I)^{-1}\nonumber\\
&=&2(W(\la)+I)^{-1}.\label{GW}
\end{eqnarray}
Therefore, $F$ is given by $F(\la)=\Theta(\la)(I-\la G(\la))$,
$\la\in\BD$. In particular, $F$ is analytic on $\BD$ and, since
$G\in\eS(\sE,\sE)$, we obtain that $\Theta(\la)=F(\la)(I-\la
G(\la))^{-1}$ for all $\la\in\BD$. Then the definition of $\Theta$
shows that (\ref{Garep}) is satisfied. Since $G$ is the inverse
Cayley transform of $W$, the function $W$ must be the Cayley
transform of $G$. Hence, using (\ref{eq:realpart}) with $G$ in
place of $C$, the real part of $W$ is given by
\[
\Re W(\la)=(I-\la G(\la))^{-*}(I-|\la|^2G(\la)^*G(\la))(I-\la
G(\la))^{-1},\quad \la\in\BD.
\]
Then for each $\la\in\BD$ we have
\[
\left(I - \lambda G(\lambda)\right)^{-*}F(\lambda)^*F(\lambda)
\left(I - \lambda G(\lambda)\right)^{-1}
=
\Theta(\lambda)^*\Theta(\lambda)  \leq \Re W(\lambda)
\]
\[
= \left(I - \lambda G(\lambda)\right)^{-*}
\left(I - |\lambda|^2 G(\lambda)^*  G(\lambda)\right)
\left(I - \lambda G(\lambda)\right)^{-1}.
\]
Thus $F(\lambda)^*F(\lambda) + |\lambda|^2  G(\lambda)^*G(\lambda) \leq I$
for all $\lambda\in\mathbb{D}$. In other words, $\col[F,\la G]$ is in
$\eS(\sE,\sY\oplus\sE)$.
Using the maximum principle for analytic functions from $\sE$ to
$\sY \oplus \sE$ we see that $\col[F,G]$ is in $\eS(\sE,\sY\oplus\sE)$.

Note that $C$ and $W$ uniquely determine each other, by Theorem
\ref{param}, and $W$ and $G$ determine each other uniquely because
$G$ is the inverse Cayley transform of $W$. Hence $C$ and $G$
determine each other uniquely. In other words, the map $J_\Gamma$
is one to one.

To prove the surjectivity, let us assume that $\col[F,G]$ is in
$\eS(\sE,\sY\oplus\sE)$ and satisfies (\ref{Garep}). Then $G$ is a
function  in $\eS(\sE,\sE)$. Let $W$ be the Cayley transform of
$G$. Then $W$ is positive real and $W(0) = I$. Moreover, for each
$\lambda$  in $\mathbb{D}$ we have
\begin{eqnarray*}
\Theta(\lambda)^*\Theta(\lambda)
&=& \left(I - \lambda G(\lambda)\right)^{-*}F(\lambda)^*F(\lambda)
\left(I - \lambda G(\lambda)\right)^{-1}\\
&\leq&  \left(I - \lambda G(\lambda)\right)^{-*}
\left(I -  G(\lambda)^*  G(\lambda)\right)
\left(I - \lambda G(\lambda)\right)^{-1} \\
&\leq&  \left(I - \lambda G(\lambda)\right)^{-*}
\left(I - |\lambda|^2 G(\lambda)^*  G(\lambda)\right)
\left(I - \lambda G(\lambda)\right)^{-1}
 = \Re W(\lambda).
\end{eqnarray*}
Thus $W$ is a $\eL(\sE,\sE)$-valued positive real function  that
satisfies (\ref{eq:bound1}), and we can apply Theorem \ref{param}
to show that $W$ is given by (\ref{defW2}) for some function $C$
in $\eS(\sD_\Gamma,\sD_\Gamma)$. Since $W$ is the Cayley
transform of $G$, we have
\[
G(\la)=\la^{-1}(W(\la)-I)(W(\la)+I)^{-1},\quad\la\in\BD.
\]
Furthermore, (\ref{Garep}) and (\ref{GW}) yield
$F(\la)=2\Theta(\la)(W(\la)+I)^{-1}$ for all $\la\in\BD$.
We see that $\col[F,G]$ is equal to
$J_\Gamma C$.

The final statement about uniqueness is trivial, because $\Gamma$
is an isometry if and only if $D_\Gamma$ is a zero operator. \epr

\medskip Note that for the case when $\sE=\sD_A$ and $\sY=\sD_{T'}$, the map
$J_\Gamma$ in Theorem \ref{th:H2} is precisely the map $J_\Gamma$
in (\ref{defJ}).

Next, in order to deal with the constraint in (\ref{altfundeq2})
and to prove the main theorems, we first prove the following
result.

\begin{proposition}\label{pr:const}
Consider the data set $\{A,T^\prime,U^\prime,R,Q\}$ with
$U^\prime$ being given by $(\ref{szns})$. Let $\Gamma$ be a
contraction from $\sD_A$ into $H^2(\sD_{T'})$, and let $C$ be a
function in $\eS(\sD_\Gamma,\sD_\Gamma)$. Define functions $F$ and
$G$ by $\col[F,G] = J_\Gamma C$ using $(\ref{defJ})$ and
$(\ref{defW})$. Then $\{F,G\}$ is a Schur pair associated with the
given data set if and only if $\Gamma$ satisfies
$(\ref{altfundeq2})$ and $C$ belongs to
$\eS_\Omega(\sD_\Gamma,\sD_\Gamma)$.
\end{proposition}

\bpr
Let $\Theta$ be the symbol of $\Gamma$, that is, $\Theta(\la)d=(\Gamma d)(\la)$
for all $d\in\sD_A$ and all $\la\in\BD$.
Observe that $W$ in (\ref{defW}) can be rewritten as
\[
W(\la)=\Gamma^*(I-\la S^*)^{-1}(I+\la S^*)\Gamma +D_\Gamma(I-\la
C(\la))^{-1}(I+\la C(\la))D_\Gamma,\quad \la\in\BD.
\]
Since $\Gamma^*\Gamma+D_\Gamma^2=I$, we obtain
\begin{eqnarray}
W(\la)-I&=&2\la\Gamma^*(I-\la S^*)^{-1}S^*\Gamma +2\la
D_\Gamma(I-\la C(\la))^{-1}C(\la)D_\Gamma,\quad \la\in\BD,\label{Weq1}\\
\noalign{\vskip4pt} W(\la)+I&=&2\Gamma^*(I-\la S^*)^{-1}\Gamma
+2D_\Gamma(I-\la C(\la))^{-1}D_\Gamma,\quad \la\in\BD.\label{Weq2}
\end{eqnarray}
We divide the remaining part of the proof into two parts.

\smallskip\noindent\textit{Part 1}. First, assuming that $\Gamma$
satisfies (\ref{altfundeq2}), we show that $G(\la)|\sF=\om_2$ for
all $\la\in\BD$ if and only if $C$ belongs to
$\eS_\Om(\sD_\Gamma,\sD_\Gamma)$. So assume that $\Gamma$
satisfies (\ref{altfundeq2}). Using (\ref{Weq1}) and (\ref{Weq2})
we see that for  $f\in\sF$ and $\la\in\BD$ we have
\begin{eqnarray*}
\la^{-1}(W(\la)-I)f
&=& 2\Gamma^*(I-\la S^*)^{-1}S^*\Gamma f
+2D_\Gamma(I-\la C(\la))^{-1}C(\la)D_\Gamma f\\
&=& 2\Gamma^*(I-\la S^*)^{-1}S^*(E\om_1f+S\Gamma\om_2f)
+2D_\Gamma(I-\la C(\la))^{-1}C(\la)D_\Gamma f\\
&=& 2\Gamma^*(I-\la S^*)^{-1}\Gamma\om_2f
+2D_\Gamma(I-\la C(\la))^{-1}C(\la)D_\Gamma f\\
&=& (W(\la)+I)\om_2f-2D_\Gamma(I-\la C(\la))^{-1}D_\Gamma\om_2f
+2D_\Gamma(I-\la C(\la))^{-1}C(\la)D_\Gamma f\\
&=& (W(\la)+I)\om_2f+2D_\Gamma(I-\la C(\la))^{-1}
(C(\la)D_\Gamma f-D_\Gamma\om_2f).
\end{eqnarray*}
Since $G$ is defined as the inverse Cayley transform of $W$, we
obtain for $f\in\sF$ and   $\la\in\BD$ that
\begin{eqnarray}
G(\la)f
&=&\la^{-1}(W(\la)-I)(W(\la)+I)^{-1}f
=(W(\la)+I)^{-1}\la^{-1}(W(\la)-I)f \label{G|F}   \\
&=&\om_2f+2(W(\la)+I)^{-1}D_\Gamma(I-\la C(\la))^{-1}
(C(\la)D_\Gamma f-D_\Gamma\om_2f).\nonumber
\end{eqnarray}
If, in addition,  $C\in\eS_\Omega(\sD_\Gamma,\sD_\Gamma)$, then
$C(\la)D_\Gamma f=D_\Gamma \om_2 f$ for all $f\in\sF$ and all
$\la\in\BD$. In this case, the last term in (\ref{G|F}) vanishes.
In other words, $G(\la)|\sF = \om_2$ for all $\la\in\BD$.

Conversely, if $G(\la)|\sF=\om_2$ for all $\la\in\BD$, then
(\ref{G|F}) shows that
\[
(W(\la)+I)^{-1}D_\Gamma(I-\la C(\la))^{-1} (C(\la)D_\Gamma
f-D_\Gamma\om_2f)=0, \quad f\in\sF, \la\in\BD.
\]
Since $I-\la C(\la)$ is an invertible operator on $\sD_\Gamma$ for
all $\la\in\BD$ and $D_\Gamma|\sD_\Gamma$ is one to one, this
implies that $C(\la)D_\Gamma|\sF-D_\Gamma\om_2=0$ for all
$\la\in\BD$. Therefore $C$ is in
$\eS_\Omega(\sD_\Gamma,\sD_\Gamma)$. This verifies our claim.

\smallskip\noindent\textit{Part 2}. In this part we prove our proposition.
First assume that $\Gamma$ satisfies (\ref{altfundeq2}) and
$C\in\eS_\Omega(\sD_\Gamma,\sD_\Gamma)$. The result of the first
part  shows that $G(\la)|\sF=\om_2$ for all $\la\in\BD$. Since
$\col[F,G]=J_\Gamma C$, and $I-\lambda
G(\lambda)=2(I+W(\lambda))^{-1}$, we have $F(\la)=\Theta(\la)(I-\la
G(\la))$ for all $\la\in\BD$. Thus
\begin{eqnarray*}
 F(\la)f
&=&\Theta(\la)f-\la\Theta(\la)G(\la)f =(\Gamma
f)(\la)-\la\Theta(\la)\om_2f \\
\noalign{\vskip4pt} &=&\om_1f+\la(\Gamma \om_2f)(\la)-\la(\Gamma
\om_2f)(\la) =\om_1f, \quad f\in\sF, \la\in\BD.
\end{eqnarray*}
This proves that $\{F,G\}$ is a Schur pair.

Conversely,  assume that $\{F,G\}$ is a Schur pair associated with
the given data set. Since $F(\la)|\sF=\om_1$ and
$F(\la)=\Theta(\la)(I-\la G(\la))$ for all $\la\in\BD$, we obtain
for all $f\in\sF$ and all $\la\in\BD$ that
\[
\om_1 f
=F(\la)f=\Theta(\la)f-\la\Theta(\la)G(\la)f
=(\Gamma f)(\la)-\la\Theta(\la)\om_2f
=(\Gamma f)(\la)-\la(\Gamma\om_2f)(\la).
\]
In other words,  $\Gamma$ satisfies the constraint in
(\ref{altfundeq2}). Using this along with $G(\la)|\sF=\om_2$ for
all $\la\in\BD$, the result of the first part shows  that $C$ is
in $\eS_\Omega(\sD_\Gamma,\sD_\Gamma)$. \epr

\smallskip\noindent
{\bf Proof of Theorem  \ref{thrc}.} Let $\{F,G\}$ be a Schur pair
associated with the given data set. Then Theorem \ref{th:H2} and
Proposition \ref{pr:const}  show that $\Gamma$ given by
(\ref{sols2}) is a contraction from $\sD_A$ into $H^2(\sD_{T'})$
satisfying (\ref{altfundeq}). Hence $B$ given by (\ref{sols1}) is
a solution to the RCL problem.

Conversely, assume that $B$ is a solution to the RCL problem. Then
$B$ admits a matrix representation of the form (\ref{sols1}),
where $\Gamma$ is a contraction from $\sD_A$ into $H^2(\sD_{T'})$
satisfying (\ref{altfundeq}). Recall that the set
$\eS_\Omega(\sD_\Gamma,\sD_\Gamma)$ is not empty. Let $C$ be any
function in $\eS_\Omega(\sD_\Gamma,\sD_\Gamma)$.  Then we obtain
from Proposition \ref{pr:const} that the pair of functions
$\{F,G\}$ given by $\col[F,G]=J_\Gamma C$  form a Schur pair
associated with the given data set. Moreover, Theorem \ref{th:H2}
shows that $\Gamma$ satisfies (\ref{sols2}). \epr

\smallskip\noindent
{\bf Proof  Theorem  \ref{th2}.}
Assume  that  $B$ is a solution to the RCL problem. Recall that  $B$ admits a
matrix representation of the form (\ref{sols1}), where $\Gamma$ is a
contraction from $\sD_A$ into $H^2(\sD_{T'})$ satisfying the constraint in (\ref{altfundeq}).
Then Proposition \ref{pr:const} implies that $J_\Gamma$ maps
$\eS_\Omega(\sD_\Gamma,\sD_\Gamma)$ onto the set of Schur pairs $\{F,G\}$
such that (\ref{sols2}) holds. According
to Theorem \ref{th:H2} the map $J_\Gamma$ is one to one.
\epr

\medskip
As one may expect from the proof of Theorem \ref{CLP}, under
appropriate additional conditions on the data set
$\{A,T^\prime,U^\prime,R,Q\}$, the formula describing all
solutions in Theorem \ref{thrc} will yield a proper
parametrization, that is, the relation between the Schur pair
$\{F,G\}$ and the solution $B$ is one to one. We plan to come back
to this question and the related question of uniqueness of the
solution in a future publication.


\noindent
A.E. Frazho, Department of Aeronautics and Astronautics,
Purdue University,
West Lafayette, IN 47907, USA,
\emph{ e-mail:} frazho@ecn.purdue.edu.

\smallskip\noindent
S. ter Horst (corresponding author), Afdeling Wiskunde,
Faculteit der Exacte Wetenschappen,
{Vrije Universiteit,
De Boelelaan 1081a, 1081 HV Amsterdam, The Netherlands,
\emph{ e-mail:}}\linebreak
S.ter.Horst@few.vu.nl.

\smallskip\noindent
M.A. Kaashoek, Afdeling Wiskunde,
Faculteit der Exacte Wetenschappen,
Vrije Universiteit,
De Boelelaan 1081a, 1081 HV Amsterdam, The Netherlands,
\emph{ e-mail:} ma.kaashoek@few.vu.nl.

\end{document}